\documentclass[12pt]{article}
\usepackage{amssymb}
\usepackage{amsfonts}
\usepackage{amsmath}
\usepackage[usenames]{color}
\usepackage{mathrsfs}
\usepackage{amsfonts}
\usepackage{amssymb,amsmath}
\usepackage{CJK}
\usepackage{cite}
\usepackage{cases}
\usepackage{amsthm}\usepackage{multirow}
\usepackage{graphicx}
\usepackage{float}
\usepackage{makecell}
\usepackage{latexsym}
\usepackage{CJK}

\def\cl{\centerline}

\def\vs{\vspace*}

\def\ni{\noindent}

\numberwithin{equation}{section}
\newtheorem{theo}{Theorem}[section]
\newtheorem{defi}[theo]{Definition}
\newtheorem{coro}[theo]{Corollary}
\newtheorem{lemm}[theo]{Lemma}
\newtheorem{exam}[theo]{Example}

\newtheorem{rema}[theo]{Remark}

\newtheorem*{thmA}{Theorem A}

\begin{document}
\begin{center}
\cl{\large\bf \vs{8pt} Isoparametric functions on Finsler space forms\,}
\footnote {$^*\,$ Project supported by NNSFC (No.12401056) and NNSFC (No. 12261034).
\\\indent\ \ $^\dag\,$ chenyl@ahnu.edu.cn

}
\cl{Yali Chen$^1$$^\dag\,$, Qun He$^2$}

\cl{\small 1 School of Mathematics and Statistics, Anhui Normal University, Wuhu, Anhui, 241000, China.}
\cl{\small
2 School of Mathematical Sciences, Tongji University, Shanghai,
200092, China.}
\end{center}

{\small
\parskip .005 truein
\baselineskip 3pt \lineskip 3pt

\noindent{{\bf Abstract:}
In this paper, we prove that each connected regular level hypersurface of a transnormal function on a Finsler space form $(N(c),F)$ has constant principal curvatures under suitable conditions. This result extends Miyaoka's refinement of Q. M. Wang's classical theorem on transnormal functions in Riemannian space forms to the Finsler setting. On forward Finsler space forms $(\overrightarrow{N}(c),F)$, we discuss the relationship between umbilic hypersurfaces and isoparametric hypersurfaces, and construct global and local isoparametric functions using the distance function.
\vs{5pt}

\ni{\bf Key words:}
Finsler space form; Transnormal function; Isoparametric function; Umbilic; Distance function.}

\ni{\it Mathematics Subject Classification (2020):} 53C60, 53C42, 53B40.}
\parskip .001 truein\baselineskip 6pt \lineskip 6pt
\section{Introduction}
~~~~The study of isoparametric theory has a long history. E. Cartan initiated the systematic study of isoparametric hypersurfaces in real space forms of constant sectional curvature~\cite{C}. Subsequently, many mathematicians have contributed to the theory of isoparametric functions and the classification of isoparametric hypersurfaces in real space forms~\cite{C1,TP,CS1}. The fundamental structure result for a general Riemannian manifold $N$ was established by Q. M. Wang in~\cite{QM}. It shows that the focal varieties of a transnormal function $f$ are smooth submanifolds of $N$. Furthermore, each regular level set of $f$ is a tube over either of the focal varieties. In particular, Q. M. Wang obtained the following remarkable result:
\begin{thmA}~\cite{QM}
Transnormal functions on $\mathbb{S}^{n}$ and $\mathbb{E}^{n}$ are isoparametric.
\end{thmA}

However, no proof was provided in \cite{QM}. Miyaoka \cite{M} subsequently revisited this result and pointed out that a regular level set of a transnormal function need not be connected. Consequently, even if each connected component has constant mean curvature, the mean curvature may differ between components of the same level set, and hence the transnormal function itself need not be isoparametric. Miyaoka refined Wang's theorem by proving that every connected regular level hypersurface of a transnormal function on $\mathbb{S}^n$ or $\mathbb{E}^n$ is an isoparametric hypersurface, without necessarily requiring the function itself to be isoparametric. Moreover, she extended this result to hyperbolic spaces under an additional condition on the principal curvatures and constructed a counterexample showing that such a condition cannot be omitted in general. Motivated by Wang's fundamental work and Miyaoka's refinement, we investigate the corresponding problem on Finsler spaces.

In Finsler geometry, the concepts of isoparametric and transnormal functions in Finsler geometry were initially introduced in \cite{HYS}. Let $(N,F)$ be an
$n$-dimensional Finsler manifold. A non-constant $C^{1}$ function $f$
on $N$ is called \textit{isoparametric} if it is smooth on $N_{f}=\{x\in N\mid df(x)\neq0\}$ and there exist a smooth function $a(t)$ and a continuous function $b(t)$ defined on $f(N_{f})$ such that
\begin{equation}\label{1.1} \left\{\begin{aligned}
&F^{2}(\nabla f)=a(f),\\
&\Delta f=b(f),
\end{aligned}\right.
\end{equation}
where $\nabla f$ denotes the gradient of $f$, which is defined by means of the Legendre transformation, and $\Delta f$ is a nonlinear Finsler-Laplacian of $f$. Each regular level hypersurface $M_{t}=f^{-1}(t)$ is called an \textit{isoparametric hypersurface}. The function satisfying the first equation in (\ref{1.1}) is called a \textit{transnormal function}. In fact, there are two definitions of Laplacian, one depending on the volume form and one not. Their choice refines the isoparametric concept into isoparametric and $d\mu$-isoparametric\cite{HCY}.

In analogy with Riemannian geometry, a complete, simply connected Finsler manifold of constant flag curvature $c$ is called a Finsler space form and denoted by $N(c)$. If it is only forward complete, we call it a forward Finsler space form and denote it by $\overrightarrow{N}(c)$ for simplicity. Research on isoparametric hypersurfaces in Finsler space forms includes classifications for several specific types, as well as studies on the Cartan identity in general Finsler space forms~\cite{HYS1,HDY,HD,HCY}.

The fundamental structure result in \cite{QM} was generalized to Finsler geometry in~\cite{CH}. Building upon the results established in \cite{CH}, in this paper, we aim to generalize Theorem A to Finsler geometry. In Finslerian case, the non-reversibility of the metric makes geodesics irreversible. Consequently, an $f$-segment, upon hitting the focal variety, ceases to be an $f$-segment with direction reversed. Moreover, the $f$-segments naturally induce a locally Riemannian metric on the manifold, which is undefined on the focal varieties. Additionally, the functions $f$, $a(t)$ and $b(t)$ in (\ref{1.1}) are smooth only at $N_{f}$ and $f(N_{f})$, and therefore exhibit lower regularity than compared to their Riemannian counterparts. These features pose significant challenges to the study of global properties of transnormal functions in Finsler geometry.

For a connected Finsler manifold $(N,F)$, denote the focal variety of $f$ by $M_{\pm}$, where $M_{+}$ (resp. $M_{-}$) is the set that $f$ attains its global maximum (resp. minimum) value (if either of them exists). Let $f$ be a $C^{2}$ transnormal function satisfying $F^{2}(\nabla f)=a(f)$, where $a(t)\neq 0$ for $t\in\operatorname{int}f(N)$ and $a'(t)\neq 0$ for $t\in\partial f(N)\cap f(N)$~(in the one-sided sense). We call such a function $f$ an \textit{appropriate transnormal function}. Then we have

\begin{theo}\label{thm3}
If $f$ is an appropriate transnormal function and $-f$ is a transnormal function on $(N(c),F)$, then every connected regular level hypersurface of $f$ has constant principal curvatures when $c\geq0$,~or when $c<0$ and every principal curvature $\kappa$ of every connected regular level hypersurface satisfies $|\kappa|\geq\sqrt{-c}$. In particular, if every regular level set is connected, then $f$ is isoparametric under the same curvature conditions. If additionally the $\textbf{S}$-curvature is constant, the same conclusion holds with isoparametric replaced by $d\mu$-isoparametric.
\end{theo}

\begin{rema}
Compared with Theorem A, there are more conditions on $f$ and $-f$ in Theorem 1.1 which hold in Riemannian case naturally. The first conclusion of Theorem \ref{thm3} is leafwise: it concerns each connected component of a regular level set. If $M_{\pm}$ are submanifolds and $f$ has no critical value in $int f(N)$, then the condition in Theorem \ref{thm3} that $f$ is an appropriate transnormal function can be weakened to just requiring $f$ to be a transnormal function. Additionally, in Finsler geometry, the two choices of nonlinear Laplacian (with or without volume form) give rise to $d\mu$-isoparametric and isoparametric notions, respectively (See Subsection 2.1 for details). There exist functions which satisfy Theorem 1.1(See Example \ref{exam}).
\end{rema}

\begin{theo} \label{thm03}
Let $f$ be a transnormal function on a forward Finsler space form $(\overrightarrow{N}(c),F)$ with $\dim N=n\geq3$. Suppose that $f$ has no critical values in $\operatorname{int}f(N)$. If some regular level hypersurface $M_{0}=f^{-1}(t_{0})$ is connected and totally umbilic, then every regular level hypersurface of $f$ is connected and totally umbilic with constant principal curvature. Consequently, $f$ is isoparametric.
\end{theo}

\begin{theo} \label{thm02}
Let $r(x)=d_{F}(p,x)$ be a distance function from any $p\in N$ on $(\overrightarrow{N}(c),F)$. If $c>0$, then $f = \frac{1}{2}r^2(x)$ is a local isoparametric function on the cut-domain of $p$, with totally umbilical local isoparametric hypersurfaces. In particular, for the case of a standard Finsler sphere, $f = -\cos(\sqrt{c}r)$ serves as a global isoparametric function. If $c\leq0$, then $f = \frac{1}{2}r^2(x)$ constitutes a global isoparametric function on $(\overrightarrow{N}(c),F)$.
\end{theo}

\begin{rema}
If $(\overrightarrow{N}(c),F)$ additionally has constant $\textbf{S}$-curvature, then the isoparametric functions in Theorems 1.3 and 1.4 are actually $d\mu$-isoparametric.
\end{rema}
\noindent
\section{Preliminaries}
~~~~In this section, we will give some definitions and lemmas that will be used in the proof of our main results. Here and from now on, we will use the following convention of index ranges unless other stated,
$$1\leq i,j,\dots\leq n,\ \ \ \ 1\leq a,b,\dots\leq n-1.$$
\subsection{Finsler manifolds}
~~~~Let $(N,F)$ be an $n$-dimensional manifold and let $TN=\cup_{x\in N}T_xN$ be the tangent bundle over $N$ with local coordinates $(x,y)$, where $x=(x^i)$ and $y=y^i\frac{\partial}{\partial x^{i}}$. The fundamental form~$g$ of~$(N,F)$ is given by
\begin{equation*}
g=g_{ij}(x,y)dx^{i} \otimes dx^{j}, ~~~~~~~g_{ij}(x,y)=\frac{1}{2}[F^{2}] _{y^{i}y^{j}}.
\end{equation*}

The geodesics $\gamma=\gamma(t)$ of $F$ in local coordinates $(x^i)$ are characterized by
\begin{align}\label{2.1}
\frac{\textmd{d}^2x^i}{\textmd{d}t^2}+2G^i(x,\frac{\textmd{d}x}{\textmd{d}t})=0,
\end{align}
where $(x^i(t))$ are the coordinates of $\gamma(t)$ and $G^i=G^i(x,y)$ are defined by
$$G^i=\frac{g^{il}}{4}\{[F^2]_{x^ky^l}y^k-[F^2]_{x^l}\},$$
which are called the \textit{geodesic coefficients}.

Let $\pi: TN\rightarrow N$ be the natural projection. There exists the unique
\emph{Chern connection}~$\nabla$ on the pull-back bundle~$\pi^{\ast}TN$. Denote $\frac{\delta}{\delta x^{i}}:=\frac{\partial}{\partial x^{i}}-N^{j}_{i}\frac{\partial}{\partial y^{j}}$. The connection satisfies $\nabla
\frac{\partial}{\partial x^i}=\Gamma^{i}_{jk}dx^k\otimes\frac{\partial}{\partial x^{j}}$ , where the coefficients are given by
$$\Gamma^{i}_{jk}=\frac{1}{2}g^{im}(\frac{\delta g_{jm}}{\delta x^{k}}+\frac{\delta g_{km}}{\delta x^{j}}-\frac{\delta g_{jk}}{\delta x^{m}}).$$
The Riemannian curvature tensor is defined by
$$R^{i}_{k}:=y^{j}R^{i}_{j~kl}y^{l},\ \ \ \ R_{jk}:=g_{ij}R^{i}_{k}=y^{i}R_{ijkl}y^{l},$$
where
\begin{align}\label{2.6.6.6}
R^{i}_{j~kl}=\frac{\delta\Gamma^{i}_{jl}}{\delta x^{k}}-\frac{\delta\Gamma^{i}_{jk}}{\delta x^{l}}+\Gamma^{i}_{km}\Gamma^{m}_{jl}-\Gamma^{i}_{lm}\Gamma^{m}_{jk}.
\end{align}
The non-Riemannian curvature tensor is defined by
\begin{align}\label{2.6}
P^{i}_{j~kl}=-\frac{\partial\Gamma^{i}_{jk}}{\partial y^{l}}.
\end{align}
Denote
\begin{align}\label{2.6.6}
L^{i}_{~kl}=-y^{j}P^{i}_{j~kl}=(G^{i})_{y^{k}y^{l}}-\Gamma^{i}_{kl}.
\end{align}
For a vector $v=v^{i}\frac{\partial}{\partial x^{i}}$ satisfying $g_{ij}v^{i}v^{j}=1$ and $g_{ij}y^{i}v^{j}=0$, the flag curvature of $(N,F)$ is defined as $$K(y,v):=F^{-2}R_{jk}v^{j}v^{k}.$$

For $X=X^{i}\frac{\partial}{\partial x^{i}}\in\Gamma(TN)$, the covariant derivative of $X$ along $v=v^{i}\frac{\partial}{\partial x^{i}}\in T_{x}N$ with respect to a reference vector $w\in T_{x}N\setminus\{0\}$ is defined by
\begin{align}\label{new}
\nabla^{w}_{v}X(x)=\{v^{j}\frac{\partial X^{i}}{\partial x^{j}}(x)+\Gamma^{i}_{jk}(w)v^{j}X^{k}(x)\}\frac{\partial}{\partial x^{i}}.
\end{align}
Consequently, the geodesic equation (\ref{2.1}) can be expressed as $\nabla^{\dot{\gamma}}_{\dot{\gamma}}\dot{\gamma}\equiv0$. If a vector field $J(s)$ along $\gamma(s)$ satisfies
\begin{align}\label{0.0}
\nabla_{\dot{\gamma}}^{\dot{\gamma}}\nabla_{\dot{\gamma}}^{\dot{\gamma}}J+\textbf{R}_{\dot{\gamma}}(J)=0,
\end{align}
then $J(s)$ is called the \textit{Jacobi field} along $\gamma(s)$, where $\textbf{R}_{y}:=R^{i}_{k}(y)\frac{\partial}{\partial x^{i}}\otimes dx^{k}$. Let $(N,F)$ be a forward complete Finsler space. For $y\in S_{x}N=\{y\in TN\setminus\{0\}\big|F(y)=1\}$, define $\textbf{c}_{y}>0$ to be the first number $d>0$ such that there exists a Jacobi field $J(s)$ along $\gamma(s)=\exp_{x}(sy)$, $0<s<d$, satisfying $J(0)= 0 = J(d)$. $\textbf{c}_{y}$ is called the \textit{conjugate value} of $y$. From the Cartan-Hadamard Theorem, the condition $K \leq 0$ implies that the conjugate value $\mathbf{c}_y = \infty$ for any $y \in SM$. Consequently, if $N$ is still simply connected, the exponential map $\exp_x: T_xN \to N$ is not only non-singular for any $x \in M$, but also constitutes a global $C^1$-diffeomorphism\cite{GTM}. Let $\textbf{i}_{y}$ be the supremum of $r>0$ such that $\gamma(s)$ is minimizing on $[0,r]$. We define the \textit{cut-domain} of $x$ as $\mathcal{D}_{x}:=\{\exp_{x}(ty),0\leq t<\textbf{i}_{y},y\in S_{x}N\}$\cite{SZ}.

Let~${\mathcal L}:TN\rightarrow T^{\ast}N$ denote the \textit{Legendre transformation}, satisfying~${\mathcal L}(\lambda
y)=\lambda {\mathcal L}(y)$ for all~$\lambda>0,~y\in TN$. Moreover,
\begin{align*}
\mathcal L^{-1}(\xi)=F^{*}(\xi)[F^{*}]_{\xi^{i}}(\xi)\frac{\partial}{\partial x^{i}},\ \ \forall \xi\in T^{*}N\setminus \{0\}, \ \ \mathcal L^{-1}(0)=0,
\end{align*}
where $F^{*}$ is the dual metric of $F$. For a smooth function~$f: N\rightarrow \mathbb{R}$, the \textit{gradient vector field} of~$f$ at~$x$ is defined as~$\nabla f(x)={\mathcal
L}^{-1}(df(x))\in T_{x}N$. The Hessian of $f$ is defined by
\begin{align}\label{2.3}
\textrm{Hes}f(X,Y):=g_{\nabla f}(\nabla^{2}f(X),Y)=X(Yf)-(\nabla^{\nabla f}_{X}Y)f,
\end{align}
where $X,Y\in TN_{f}$, $\nabla^{2}f(X)=\nabla^{\nabla f}_{X}\nabla f\in TN_{f}$.
The \emph{Laplacian} of~$f$ can be defined
by
\begin{equation}\label{2.21}
\hat{\Delta} f=\textmd{tr}_{g_{_{\nabla f}}}(\nabla^{2}f).
\end{equation}
Furthermore, the \emph{Laplacian} of~$f$ with respect to the volume form~$d\mu=\sigma(x)dx=\sigma(x)dx^{1}\wedge dx^{2}\wedge\cdots\wedge dx^{n}$ can be represented as
\begin{equation}\label{2.10}
\Delta_{\sigma} f=\textmd{div}_{\sigma}(\nabla f)=\frac{1}{\sigma}\frac{\partial}{\partial x^{i}}(\sigma g^{ij}(\nabla f)f_{j})=\hat{\Delta} f-\textbf{S}(\nabla f),
\end{equation}
where
$$\textbf{S}(x,y)=\frac{\partial G^{i}}{\partial y^{i}}-y^{i}\frac{\partial}{\partial x^{i}}(\ln \sigma(x))$$
is the \emph{$\mathbf{S}$-curvature}~\cite{SZ}.
\subsection{Anisotropic submanifolds}
~~~~Let $(N,F)$ be an $n$-dimensional Finsler manifold and $\phi: M\to(N, F)$ be an $(n-1)$-dimensional immersion.
Let $(x,y)\in TN$ have local coordinates $(x^{i},y^{i})$ and $(u,v)\in TM$ have local coordinates $(u^{a},v^{a})$. We have the following relations:
\begin{align}\label{2.X}
x^{i}=\phi^{i}(u),\ \ \ \ y^{i}=\phi^{i}_{a}v^{a},\ \ \ \ \phi^{i}_{a}=\frac{\partial\phi^{i}}{\partial u^{a}}.
\end{align}
For simplicity, we identify $\phi(x)$ with $x$ and $d\phi X$ with $X$. Define
$$\mathcal{V}(M)=\{(x,\xi)~|~x\in M,\xi\in T_x^{*}N,\xi (X)=0,\forall X\in T_xM\},$$
which is called the
\textit{normal bundle} of $M$ \cite{SZ1}. Let $\tau\in\mathcal{N}M={\mathcal L}^{-1}(\mathcal{V}(M))\subset TN$ be a vector satisfying $g_{\tau}(\tau,X)=0$ for any $X\in TM$. Here, $\tau$ is called a \textit{normal vector field} of $M$, and the locally Riemannian manifold $(M,\hat{g}=\phi^{*}g_{\tau})$ is called an \textit{anisotropic submanifold} of $(N,F)$. Let $\textbf{n}\in\mathcal{N}^0M={\mathcal L}^{-1}(\mathcal{V}^0(M))$ be the unit normal vector field, where $\mathcal{V}^0(M)=\{\nu\in \mathcal{V}(M)|F^*(\nu)=1\}$ denotes the unit normal bundle. In general, $\mathcal L^{-1}(-\nu)\neq-\mathcal L^{-1}(\nu)$. The shape operator ~${A}_{\textbf{n}}:T_xM\rightarrow T_xM$ is defined by
\begin{align}\label{2.4}
{A}_{\textbf{n}}X=-(\nabla^{\textbf{n}}_{X}\textbf{n})^{\top}.
\end{align}
The eigenvalues of ${A}_{\textbf{n}}$, denoted by $k_{1}, k_{2}, \cdots, k_{m}$, are called \emph{the principal curvatures} of $M$ with respect to $\textbf{n}$, and $\hat{H}_{\textbf{n}}=\sum\limits_{a=1}^{m}k_{a}$ is the \textit{anisotropic mean curvature} of $M$. With $\textbf{n}\in\mathcal{N}^{0}M$ and an arbitrary volume form $d\mu$ on $(N,F)$, the $d\mu$-mean curvature $H_{\textbf{n}}$ is defined, which satisfies $H_{\textbf{n}}=\hat{H}_{\textbf{n}}+\textbf{S}(\textbf{n})$~\cite{HYS}.

If $M$ is a coorientable hypersurface, we may choose a global unit normal vector field $\mathbf n$. Consequently, $(M,\hat{g})$ is an $\textit{anisotropic hypersurface}$, and the Weingarten formula (\ref{2.4}) simplifies to
\begin{align}\label{2.5}
{A}_{\textbf{n}}X=-\nabla^{\textbf{n}}_{X}\textbf{n}.
\end{align}
\subsection{Isoparametric functions and isoparametric hypersurfaces}
~~~~More precisely, the function $f$ is called \emph{isoparametric} (resp. \emph{$d\mu$-isoparametric}) on $(N, F, d\mu)$ if (\ref{1.1}) holds with $\Delta f=\hat{\Delta} f$ (resp. $\Delta f=\Delta_{\sigma}f$), as defined in (\ref{2.21}) (resp. (\ref{2.10})). If $(N,F,d\mu)$ has constant $\mathbf{S}$-curvature, then $f$ is isoparametric if and only if it is $d\mu$-isoparametric. Geometrically, isoparametric (resp. $d\mu$-isoparametric) functions and hypersurfaces can be characterized via anisotropic mean curvatures (resp. $d\mu$-mean curvatures) and principal curvatures.
\begin{defi}\cite{HYS}
A transnormal function $f$ is isoparametric (resp. $d\mu$-isoparametric) on $(N,F,d\mu)$ if and only if each regular level hypersurface of $f$ has constant anisotropic mean curvatures (resp. $d\mu$-mean curvatures) with respect to $\textbf{n}=\frac{\nabla f}{F(\nabla f)}$.
\end{defi}

\begin{defi}\cite{HYS}\label{defi1}
Let $(N,F)$ be a Finsler manifold with constant flag curvature. Then a transnormal function $f$ is isoparametric if and only if each regular level hypersurface of $f$ has constant principal curvatures. If the $\textbf{S}$-curvature is constant, the above equivalence holds with isoparametric replaced by $d\mu$-isoparametric.
\end{defi}

In analogy with the Riemannian case, we can also define local ($d\mu$-)isoparametric functions and hypersurfaces on $(N,F)$. If for any $x\in M$, there exists a neighborhood $U$ of $x$ and an ($d\mu$-)isoparametric function $f$ defined on $U$, such that $M\cap U$ is a regular level hypersurface of $f$, then $f$ is called a local ($d\mu$-)isoparametric function and $M$ is called a local ($d\mu$-)isoparametric hypersurface\cite{HH}. Geometrically, local isoparametric (resp. $d\mu$-isoparametric) functions and hypersurfaces can also be characterized via anisotropic mean curvatures (resp. $d\mu$-mean curvatures) and principal curvatures.
\begin{defi}\cite{DC}
$M$ in a Finsler manifold is a local isoparametric (resp. $d\mu$-isoparametric) hypersurface if and only if $M$ and its nearby parallel hypersurfaces have constant anisotropic mean curvatures (resp. $d\mu$-mean curvatures) with respect to $\textbf{n}=\frac{\nabla f}{F(\nabla f)}$.
\end{defi}

\begin{defi}\cite{HH}\label{defi2}
$M$ in a Finsler manifold with constant flag curvature is a local isoparametric hypersurface if and only if its principal curvatures are all constants. If the $\textbf{S}$-curvature is constant, the above equivalence holds with isoparametric replaced by $d\mu$-isoparametric.
\end{defi}
Unlike in the Riemannian case, there exist local ($d\mu$-)isoparametric hypersurfaces that cannot be extended to global ($d\mu$-)isoparametric hypersurfaces (See Remark \ref{rema1}). However, under certain additional conditions, such extensions are guaranteed (See Remark \ref{rema2}).
\section{Transnormal functions on Finsler space forms}
\subsection{Transnormal functions and focal varieties}
~~~~For $\tau\in\mathcal{N}_{x}M$, define the normal exponential map $E:\mathcal{N}M\rightarrow N$ by $E(x,\tau)=\exp_{x}\tau$, which is $C^{\infty}$ on $TN\setminus\{0\}$. The \textit{focal points} of $M$ are defined as the critical values of the normal exponential map $E$. A focal point $p\in N$ is said to have multiplicity $m$ if $E_{*}$ at $(x,\tau)$ has nullity $m>0$. The set of focal points is called the \textit{focal set}, which is denoted by $\textmd{Foc}(f)$. The focal varieties of $f$ is denoted by $M_{\pm}$, which is the set where $f$ attains its global maximum value or global minimum value (if either of them exists). Denote $V=\{p\in N\big|\nabla f(p)=0\}$, which is the set of critical point of $f$. Indeed, $\textmd{Foc}(f)\subset M_{\pm}\subset V$.
\begin{lemm}\label{lem2}\cite{CH}
Let $f$ be a transnormal function on $(N,F)$. Suppose that $\beta$ is the only critical value of $f$ in $[\alpha, \beta]\subset f(N)$, then
\begin{equation}\label{3.2}
d(M_{\alpha}, M_{\beta})=\int^{\beta}_{\alpha}\frac{df}{\sqrt{a(f)}}
\end{equation}
and the improper integral in (\ref{3.2}) converges.
\end{lemm}
\begin{rema}
In fact, from the proof of Lemma 3.2 in \cite{CH}, we know that if $\alpha$ and $\beta$ are the only critical values of $f$ in $[\alpha, \beta]\subset f(N)$, then (\ref{3.2}) still holds.
\end{rema}
Define the reverse metric of $F$ by $\overleftarrow{F}(x,y)=F(x,-y)$, where $y\in T_{x}N$, $x\in N$. Then we have
\begin{lemm}\label{lem41}\cite{CH}
Let $(N,F)$ be a complete, connected Finsler manifold and $f$ be an appropriate transnormal function on $N$, then we have\\
(1)each focal variety of $f$ is a smooth submanifold in $N$.\\
(2)each regular level set of $f$ is a tube or a parallel hypersurface over either of the focal varieties with respect to $\overleftarrow{F}$ or $F$.
\end{lemm}
Suppose that $f$ is an appropriate transnormal function on $(N,F)$. Obviously, $V=M_{+}\cup M_{-}$. In addition, from Lemma \ref{lem41}, if $\textmd{Foc}(f)\cap M_{+}$(resp. $M_{-}$)$=\varnothing$, then $M_{+}$(resp. $M_{-}$) is a hypersurface. If $\textmd{Foc}(f)\cap M_{+}$(resp. $M_{-}$)$\neq\varnothing$, then $M_{+}$(resp. $M_{-}$) is a submanifold with codimension at least 2.

Recall that an $f$-segment is a unit-speed reparameterization of an integral curve of $\nabla f$, oriented so that $f$ is increasing. Similar to Riemannian case, we can define a \textit{maximal $f$-segment} as an $f$-segment that extends orthogonally to (and terminates at) the focal varieties. Unlike the Riemannian case, if a maximal $f$-segment in Finsler geometry orthogonally hits a focal variety, the curve starting from that intersection point is no longer an $f$-segment. This immediately raises the question: what kind of curve is it if not an $f$-segment? This question is fundamental to our study. In the following, we resolve it completely. From Lemma \ref{lem41} and the proof of Lemma 3.5 in \cite{CH}, we immediately have
\begin{lemm}\label{lem56}
Let $f$ be an appropriate transnormal function on $(N,F)$. If $\gamma(s)$ is a maximal $f$-segment satisfying $\gamma(0)\in M_{+}$ (resp. $M_{-}$), then $\gamma'(0)\in\mathcal{N}^{0}_{\gamma(0)}M_{+}$ (resp. $\mathcal{N}^{0}_{\gamma(0)}M_{-}$). Conversely, for any $\eta\in\mathcal{N}^{0}_{p}M_{+}$ (resp. $\mathcal{N}^{0}_{q}M_{-}$), there exists a maximal $f$-segment $\gamma(s)$ satisfying $\gamma(0)=p\in M_{+}$ (resp. $\gamma(0)=q\in M_{-}$) and $\gamma'(0)=\eta$.
\end{lemm}
\begin{rema}
If $f$ is only a transnormal function and $M_{\pm}$ are submanifolds, Lemma \ref{lem56} holds.
\end{rema}
Suppose that both $M_{+}$ and $M_{-}$ are nonempty. For an appropriate transnormal function $f$ satisfying $F^{2}(\nabla f)=a(f)$, from Lemma \ref{lem2}, we define
$d_{+}=d(M_{-},M_{+})=\int^{\beta}_{\alpha}\frac{dt}{\sqrt{a(t)}}.$
Let $-f$ be a transnormal function satisfying $F^{2}(\nabla(-f))=\tilde{a}(-f)$. Due to $\nabla f=0$ if and only if $\nabla (-f)=0$, we know that $-f$ has no critical values in $-f(N_{f})$. Since $M_{\pm}^{f}=M_{\mp}^{-f}$, we define $d_{-}=d(M_{+},M_{-})=\int^{\beta}_{\alpha}\frac{dt}{\sqrt{\tilde{a}(-t)}}$. From Lemma \ref{lem56} and the preceding remark, we have
\begin{lemm}\label{lem3}
Let $f$ be an appropriate transnormal function and $-f$ be a transnormal function on $(N,F)$. Suppose that both $M_+$ and $M_-$ are nonempty. Let $\gamma(s)$, $s\in[0,+\infty)$, be a unit speed geodesic such that $\gamma(s)|_{[0,r]}$ is an $f$-segment and $\gamma(r)\in M_{+}$. Then $\gamma(s)|_{[r,r+d_{-}]}$ is a maximal $(-f)$-segment with respect to $F$. Conversely, if $\gamma(s)|_{[0,r]}$ is a $(-f)$-segment and $\gamma(r)\in M_{-}$, then $\gamma(s)|_{[r,r+d_{+}]}$ is a maximal $f$-segment with respect to $F$.
\end{lemm}
\proof
Set $p=\gamma(r)\in M_{+}$. By Lemma 3.4, we have $\dot\gamma(r)\in\mathcal{N}^{0}_{p}M_{+}$ and there exists a maximal $(-f)$-segment $\sigma(s)$ satisfying
\[
\sigma(0)=p,\qquad \dot\sigma(0)=\dot\gamma(r).
\]
Both $\sigma(s)$ and $\gamma(r+s)$ are geodesics of $F$ with the same initial position and initial velocity. Hence, by uniqueness of the geodesic equation,
\[
\sigma(s)=\gamma(r+s).
\]
Its length from $M_{+}$ to $M_{-}$ is $d_{-}$. The second assertion follows in the same way by interchanging $f$ and $-f$. This completes the proof.
\endproof
\begin{rema}
If $f$ and $-f$ are transnormal functions, there are no critical levels between $M_{-}$ and $M_{+}$, and $M_{\pm}$ are submanifolds, then the conclusion of Lemma \ref{lem3} still holds.
\end{rema}
\subsection{Proof of Theorem \ref{thm3}}
\begin{lemm}\cite{HCY}\label{lem16}
Let $\phi:M  \rightarrow \overrightarrow{N}(c) $ be an immersion submanifold  and $\mathbf{n}$ be a unit normal vector to $\phi(M)$ at $x$.
Then $p=\phi_{s}(x,\mathbf{n})$ is a focal point
 of $(M ,x)$ of multiplicity $m_0>0$ if and only if there is an eigenvalue $\kappa$ of the shape operator $A_{\mathbf{n}}$ of multiplicity $m_0$ such that
\begin{align}\label{3.1}
\kappa(x)=\left\{
\begin{array}{lcl}
\frac{1}{s},     &      & {c=0,}\\
\sqrt{c}\cot \sqrt{c}s,   &      & {c>0,}\\
\sqrt{-c}\coth \sqrt{-c}s,       &      & {c<0.}
\end{array} \right.\end{align}
\end{lemm}
Let $t_{0}\in(\alpha,\beta)$ be a regular value and let $M$ be an arbitrary connected component of $f^{-1}(t_{0})$, oriented by
\[
\mathbf n_{f}=\frac{\nabla f}{\sqrt{a(f)}}.
\]
For each $x\in M$, let $\gamma_{x}:\mathbb R\rightarrow N$ be the unit speed geodesic satisfying
\[
\gamma_{x}(0)=x,\qquad \dot\gamma_{x}(0)=\mathbf n_{f}(x),
\]
and set
\[
\Phi_{s}(x)=\gamma_{x}(s),\qquad
P(M)=\{s\in\mathbb R\mid \ker(d\Phi_{s})_{x}\neq0\ \text{for some }x\in M\}.
\]
We distinguish three cases according to
the existence of $M_+$ and $M_-$. Define
$$T_M:=\{s\in\mathbb R\mid\gamma_x(s)\in M_+\cup M_-\text{ for some }x\in M\}.$$
\medskip
\noindent
\textit{Case 1. Both $M_+$ and $M_-$
are nonempty.}

Set
\[
r_{+}(t_{0})=\int_{t_{0}}^{\beta}\frac{dt}{\sqrt{a(t)}},\qquad
r_{-}(t_{0})=\int_{\alpha}^{t_{0}}\frac{dt}{\sqrt{a(t)}},\qquad
D=d_{+}+d_{-}.
\]
By Lemma \ref{lem2}, $r_{+}(t_{0})$, $r_{-}(t_{0})$ and $D$ are finite.
For $s\geq0$, the geodesic $\gamma_x$ starts from $\gamma_x(0)\in M$ and follows an $f$-segment until it reaches $M_+$ at $\gamma_x(r_+)$ and $\gamma_{x}^{'}(r_+)\in\mathcal{N}^{0}_{\gamma_{x}(r_{+})}M_{+}$. By Lemma 3.6, its continuation is a maximal $(-f)$-segment from $M_+$ to $M_-$, of length $d_-$. By the converse assertion of Lemma 3.6, it then continues as an $f$-segment from $M_-$ to $M_+$, of length $d_+$. This alternating pattern repeats along the same geodesic.
Thus, the parameters at which $\gamma_x$ reaches $M_+$ or $M_-$ for $s\geq 0$ are given by
$$T_M^+=\{r_++kD,\ r_++d_-+kD\mid k\in\mathbb N_0\}.$$

For $s<0$, let $\overleftarrow F(x,y)=F(x,-y)$ and put $\bar\gamma_x(\tau)=\gamma_x(-\tau)$. Then $\bar\gamma_x$ is a unit speed geodesic of $\overleftarrow F$. By the Legendre transformation,
$$
\nabla^{\overleftarrow F}(-f)=-\nabla^F f,
\qquad
\nabla^{\overleftarrow F}f=-\nabla^F(-f).
$$
Hence, $\bar\gamma_x$ is an $(-f)$-segment with respect to $\overleftarrow F$, starting from $x\in M$ and reaching $M_-$ at $\bar\gamma_x(r_-)$. By Lemma 3.6, it then continues alternately as an $f$-segment and a $(-f)$-segment with respect to $\overleftarrow F$, with successive lengths $d_-$ and $d_+$. Consequently,
$$
T_M^-=
\{-r_--kD,,-r_--d_--kD\mid k\in\mathbb N_0\}.
$$
Consequently,
$$T_M=T_M^+\cup T_M^-.$$
\medskip
\noindent
\textit{Case 2. Exactly one of $M_+$
and $M_-$ is nonempty.}

If $M_+\ne\varnothing$, each $\gamma_x$ follows an $f$-segment until it reaches $M_+$ at $s=r_+$, and then continues as a $(-f)$-segment. Since $M_-$ is empty, $T_M=\{r_+\}.$

If $M_-\ne\varnothing$, as in Case~1, the reversed geodesic $\bar\gamma_x$ follows a $(-f)$-segment with respect to
$\overleftarrow F$ until until $\bar\gamma_x(r_-)\in M_-$, and then continues as an $f$-segment. Since $M_+$ is empty, $T_M=\{-r_-\}.$

\medskip
\noindent
\textit{Case 3. Both $M_+$ and $M_-$
are empty.}

There are no critical values, and every $\gamma_x$ is an $f$-segment for all $s\in\mathbb R$. Hence, $T_M=\varnothing.$

In all three cases, $T_M$ is countable, and for every $x\in M$, $T_M=\{s\in\mathbb R:\gamma_x(s)\in M_+\cup M_-\}.$ Moreover, for $s\notin T_M$, the geodesic $\gamma_x$ lies in a regular $f$- or $(-f)$-segment.

Suppose that $s_0\in P(M)\setminus T_M$. Then there exist
$x_0\in M$ and a nonzero vector $X\in T_{x_0}M$ such that
$$
(d\Phi_{s_0})_{x_0}(X)=0.
$$
Choose a smooth curve $x(t)$ in $M$ satisfying
$x(0)=x_0$ and $\dot{x}(0)=X$, and consider the
geodesic variation
$$
\Psi(s,t)=\Phi_s(x(t)).
$$
Set
$$
\widetilde J(s,t)=\frac{\partial\Psi}{\partial t}(s,t),
\qquad
\widetilde T(s,t)=\frac{\partial\Psi}{\partial s}(s,t),
$$
and
$$
J(s)=\widetilde J(s,0),
\qquad
T(s)=\widetilde T(s,0).
$$
Then $J$ is a Jacobi field along $\gamma_{x_0}$,
satisfying
$$
J(0)=X\ne0,
\qquad
J(s_0)=(d\Phi_{s_0})_{x_0}(X)=0.
$$
Since $T_M$ is independent of $x\in M$, near $s_0\notin T_M$,
we have
$$
\widetilde T(s,t)=
\frac{\partial\Psi}{\partial s}(s,t)
=\mathbf{n}_f(\Psi(s,t))
$$
or
$$
\widetilde T(s,t)=
\frac{\partial\Psi}{\partial s}(s,t)
=\mathbf{n}_{-f}(\Psi(s,t)).
$$
In particular, $T(s)$ is the restriction of
this unit gradient field to $\gamma_{x_0}$. Since $J(s_0)=0$,
$$
D_T^T J(s_0)
=D_{\tilde{T}}^{\tilde{T}} \tilde{J}\mid_{t=0}
=D_{\tilde{J}}^{\tilde{T}} \tilde{T}\mid_{t=0}=D_{J(s_0)}^T \tilde{T}\mid_{t=0}=0.
$$
The uniqueness of the Jacobi equation implies
$J\equiv 0$, contradicting $J(0)=X\ne 0$.
Therefore,
$$
P(M)\subset T_M.
$$
Let
$$
\kappa_1(x)\leq\cdots\leq\kappa_{n-1}(x)
$$
be the ordered eigenvalues of the shape operator of $M$.
Since the shape operator is smooth and self-adjoint
with respect to $g_{n_f}$, each $\kappa_j$ is
continuous on $M$.

If $c>0$, by Lemma 3.8, for any $\kappa_j(x)$,
there exists some $s$ such that
$\gamma_x(s)\in Foc(f)$. Hence,
$$
\kappa_j(M)\subset
\left\{\sqrt{c}\cot(\sqrt{c}s)\mid s\in T_M\right\}.
$$

If $c=0$, by Lemma 3.8, for any $\kappa_j(x)\neq0$,
there exists some $s$ such that
$\gamma_x(s)\in Foc(f)$. Thus,
$$
\kappa_j(M)\subset
\left\{\frac{1}{s}\mid s\in T_M,\ s\neq0\right\}
\cup\{0\}.
$$

If $c<0$, by Lemma 3.8, for any
$|\kappa_j(x)|>\sqrt{-c}$, there exists some $s$
such that $\gamma_x(s)\in Foc(f)$.
Under the assumption
$|\kappa_j(x)|\geq\sqrt{-c}$, we obtain
$$
\kappa_j(M)\subset
\left\{\sqrt{-c}\coth(\sqrt{-c}s)
\mid s\in T_M\right\}
\cup\{-\sqrt{-c},\sqrt{-c}\}.
$$

Since $T_M$ is countable, each $\kappa_j(M)$ is
countable. As $M$ is connected and $\kappa_j$ is
continuous, each $\kappa_j$ must be constant on $M$.
Since the connected regular leaf $M$ was arbitrary,
every connected regular level hypersurface of $f$
has constant principal curvatures. If every regular level set is connected, Definition \ref{defi1} implies that $f$ is isoparametric. If $(N(c),F)$ additionally has constant $\textbf{S}$-curvature, the same conclusion holds for $d\mu$-isoparametric functions. This completes the proof of Theorem \ref{thm3}.
\endproof

\begin{rema}
He--Chen--Yin--Ren \cite{HCY} obtained the bound $g\leq 2$ for $c\leq 0$ under the assumption of vanishing reversible torsion, using a generalized Cartan identity. However, their focal-submanifold argument requires further justification in two respects.
First, when $c<0$, a principal curvature $\kappa_j$ gives rise to a focal point at a finite parameter only if $|\kappa_j|>\sqrt{-c}$. Second, applying the shape operator formula to every unit normal vector of a focal submanifold requires
a global normal-coverage property, which does not follow from local tube formulas alone. Under the additional transnormal assumptions on $f$ and $-f$ in Theorem \ref{thm3}, together with appropriate principal-curvature conditions, these difficulties can be addressed.
\end{rema}

If $F$ is reversible, $\tilde{a}(-f)=(F^{*}(-df))^{2}=(F^{*}(df))^{2}=F^{2}(\nabla f)=a(f)$. That is, $f$ is an appropriate transnormal function if and only if $-f$ is an appropriate transnormal function. Hence, we have
\begin{coro}
Let $f$ be an appropriate transnormal function on
$(N(c),F)$, where $F$ is reversible. Then every connected regular level hypersurface of $f$ has constant principal curvatures if $c\geq 0$, or if $c<0$ and all principal curvatures $\kappa_i$ satisfy $|\kappa_i|\geq\sqrt{-c}$. If every regular level set of $f$ is connected, then $f$ is isoparametric under the same curvature conditions. If, in addition, the $S$-curvature is constant, then $f$ is $d\mu$-isoparametric under the same curvature conditions.
\end{coro}

\begin{exam}\label{exam}
The functions constructed on standard Randers spheres in [12, Examples 3.12--3.13] satisfy the assumptions of Theorem 1.1 and are globally \((d\mu\text{-})\)isoparametric, as established in \cite{HDY}. Furthermore, Theorem 1.1 can also be applied to construct new examples of ($d\mu$-)isoparametric functions without computing their Laplacian. Let $f(x)=sin(\langle a,x\rangle)$, defined on a Minkowski space equipped with an $(\alpha,\beta)$-metric $F(y)=\alpha(y)\phi(\frac{\beta}{\alpha})$, where $\alpha$ is Euclidean and $\beta=\langle p,y\rangle$ with a constant vector $p$. Set the dual metric $F^{*}(\xi)=\alpha^{*}(\xi)\Phi(\frac{\beta^{*}}{\alpha^{*}})$, where $\alpha^{*}$ is the dual metric of $\alpha$ and $\beta^{*}$ is the dual vector field of $\beta$ with respect to $\alpha$. Assume $\langle a,a\rangle=1$ and $\langle p,a\rangle=0$. Then $F^{*}(df)=\sqrt{1-f^{2}}\Phi(0)=F^{*}(-df)$. Let $\sqrt{C}=\Phi(0)>0$. Then $a(t)=\tilde{a}(t)=C(1-t^{2})$. Thus, $\pm f$ are both appropriate transnormal functions. By Theorem \ref{thm3}, every connected regular level has constant principal curvatures. Since the regular components are parallel hyperplanes with zero anisotropic mean curvature, the component compatibility is automatic, and hence $f$ is a ($d\mu$-)isoparametric function for the BH or HT volume form.
\end{exam}

\section{Local isoparametric hypersurfaces}
\subsection{The number of distinct principal curvatures}
~~~~The conclusion of Theorem~1.1 does not necessarily hold for local transnormal functions. Indeed, by using the distance function to a hypersurface, one can construct a local appropriate transnormal function, such as $f=-\cos r$, whose regular level hypersurfaces need not have constant principal curvatures.

Moreover, the distinction between local and global results is also essential for the number $g$ of distinct principal curvatures. Ge and Ma~\cite{GM} showed that a complete isoparametric hypersurface in a Minkowski space has at most two distinct principal curvatures. However, this bound need not hold locally. We illustrate this by constructing a local isoparametric hypersurface with three distinct principal curvatures in a Minkowski space.

In \cite{HH}, the conic Finsler metric was defined and a local minimal isoparametric hypersurface was found in conic Minkowski spaces, which is not a conic Minkowski hyperplane, hypersphere, or cylinder. For a Minkowski norm equipped with the BH or HT volume, the $\textbf{\textit{S}}$-curvature vanishes. Thus, isoparametricity and $d\mu$-isoparametricity are equivalent in this setting.
\begin{lemm}\cite{HH}\label{thm17}
Let $(\mathbf{R}^3,F)$ be a 3-dimensional conic Minkowski-$(\alpha,\beta)$ space and the dual metric of $F$ be $F^{*}=\alpha^{*}\varphi(\frac{\beta^{*}}{\alpha^{*}})$, where
$$\varphi=\frac{\sqrt{b^2-(1+a^2)s^2}}{b}-\frac{as}{b}\arctan\frac{\sqrt{b^2-(1+a^2)s^2}}{as},0<|s|<\frac{b}{\sqrt{1+a^2}},$$
where $\beta^{*}=(0,0,b)$ and $a,b$ are two positive constants. Then in $(\mathbf{R}^3,F)$, helicoid $\mathbf{r}=(u\cos v,u\sin v,av)(0<u<1)$ is a local minimal isoparametric hypersurface with constant principal curvatures $\pm1$.
\end{lemm}
Let $F_{\mathrm{conic}}^*$ denote the three-dimensional
conic dual Minkowski norm in Lemma~4.1. Consider the
open cone
$$
U=\left\{\xi\in\mathbb{R}^3:
\xi_1^2+\xi_2^2>4a^2\xi_3^2,\ \xi_3>0\right\},
$$
on which $F_{\mathrm{conic}}^*$ is smooth and
strongly convex. For the helicoid
$$
\mathbf{r}(u,v)=(u\cos v,u\sin v,av),
\qquad 0<u<\frac12,
$$
its conormal directions lie in $U$.

For $n\geq 4$, denote
$\xi=(\eta,\zeta)\in\mathbb{R}^3\times\mathbb{R}^{n-3}$
and define
\begin{equation}
C_n^*:=U\times\mathbb{R}^{n-3},\qquad
F^{*2}(\eta,\zeta)
:=F_{\mathrm{conic}}^{*2}(\eta)+|\zeta|^2.
\end{equation}
Since $F_{\mathrm{conic}}^*$ is smooth and strongly
convex on $U$, $F^{*}$ defines an $n$-dimensional
conic dual Minkowski norm.
\medskip
\begin{lemm}\label{lemX}
Let $F^*$ be an $n$-dimensional conic dual Minkowski norm defined on an open cone $C_n^*\subset\mathbb{R}^n\setminus\{0\}$. Then, for every $\xi_0\in C_n^*$, there exist a regular dual Minkowski norm $\widetilde F^*$ on $\mathbb{R}^n$ and an open cone $C_0^*$ containing the ray $\mathbb{R}_{>0}\xi_0$, with $C_0^*\subset C_n^*$, such that
\[
\widetilde F^*\mid_{C_0^*}=F^*\mid_{C_0^*}.
\]
\end{lemm}
\proof
Define $F_0^{*2}(\xi)=g^{*ij}(\xi_0)\xi^i\xi^j$ on $\mathbb{R}^n$. By a direct calculation,
$$
\begin{aligned}
F_0^{*2}(\xi_0)&=F^{*2}(\xi_0),\\
(F_0^{*2})_{\xi_i}(\xi_0)&=(F^{*2})_{\xi_i}(\xi_0),\\
(F_0^{*2})_{\xi_i\xi_{j}}(\xi_0)&=(F^{*2})_{\xi_i\xi_{j}}(\xi_0)
\end{aligned}
$$
Put
$$
R(\xi)=F^{*2}(\xi)-F_0^{*2}(\xi).
$$
Obviously, we have
$$R(\xi_{0})=0,~~~~R_{\xi_{i}}(\xi_{0})=0,~~~~R_{\xi_{i}\xi_{j}}(\xi_{0})=0.$$

Fix an auxiliary Euclidean norm $|\cdot|$ on
$\mathbb{R}^n$, and let
$S^{n-1}=\{\xi\in\mathbb{R}^n:|\xi|=1\}$.
Choose a smooth function
$\rho:[0,+\infty)\to[0,1]$ such that
$\rho=1$ on $[0,1/4]$ and
$\rho=0$ on $[1,+\infty)$.
For $\varepsilon>0$, define
\begin{equation}\label{eq:angular-cutoff}
\chi_\varepsilon(\xi)
=\rho\left(
\frac{1}{\varepsilon^2}
\left|
\frac{\xi}{|\xi|}
-\frac{\xi_0}{|\xi_0|}
\right|^2
\right).
\end{equation}

Set $\widehat\xi_0=\xi_0/|\xi_0|$.
Since $C_n^*$ is open and $\xi_0\in C_n^*$,
we may choose $\varepsilon>0$ sufficiently
small such that
$$
\begin{aligned}
\operatorname{supp}
(\chi_\varepsilon|_{S^{n-1}})\subset
\left\{\xi\in S^{n-1}:
|\xi-\widehat\xi_0|\leq\varepsilon
\right\}\subset C_n^*\cap S^{n-1}.
\end{aligned}
$$
Define
$$
\widetilde R(\xi)=
\begin{cases}
\chi_\varepsilon(\xi)R(\xi),
&\xi\in C_n^*,\\
0,&\xi\notin C_n^*,
\end{cases}
$$
and
\begin{equation}\label{eq:extended-dual}
\widetilde F^{*2}(\xi)
=F_0^{*2}(\xi)+\widetilde R(\xi),
\qquad \xi\in\mathbb{R}^n\setminus\{0\}.
\end{equation}
Since $\chi_\varepsilon$ vanishes near the boundary of $C_n^*$, the function $\widetilde R$ is smooth on $\mathbb{R}^n\setminus\{0\}$.

We now consider the Hessian matrix of
$\widetilde F^{*2}$. By \eqref{eq:extended-dual},
$$
2\widetilde g^{*ij}({\xi})=2g^{*ij}(\xi_0)+\widetilde R_{\xi_i\xi_j}.
$$
We first estimate $\widetilde R_{\xi_i\xi_j}$. On $C_n^*$, by a direct calculation,
$$
\begin{aligned}
\widetilde R_{\xi_i\xi_j}
={}&\chi_\varepsilon R_{\xi_i\xi_j}
+(\chi_\varepsilon)_{\xi_i}R_{\xi_j}+(\chi_\varepsilon)_{\xi_j}R_{\xi_i}+R(\chi_\varepsilon)_{\xi_i\xi_j}.
\end{aligned}
$$

Since $R$ is two-homogeneous, we have
$$
\begin{aligned}
R(\widehat\xi_0)=\frac{1}{|\xi_0|^2}R(\xi_0)=0,~~~~R_{\xi_i}(\widehat\xi_0)=\frac{1}{|\xi_0|}R_{\xi_i}(\xi_0)=0,~~~~R_{\xi_i\xi_j}(\widehat\xi_0)
=R_{\xi_i\xi_j}(\xi_0)=0.
\end{aligned}
$$
Applying Taylor's theorem at
$\widehat\xi_0$ and using
$|\xi-\widehat\xi_0|\leq\varepsilon$
on the support of $\chi_\varepsilon$
in $S^{n-1}$, we obtain
$$
\begin{aligned}
|R(\xi)|&\leq A_0\varepsilon^3,~~~~|R_{\xi_i}(\xi)|&\leq A_1\varepsilon^2,~~~~|R_{\xi_i\xi_j}(\xi)|&\leq A_2\varepsilon,
\end{aligned}
$$
for $\xi\in
\operatorname{supp}(\chi_\varepsilon|_{S^{n-1}})$,
where $A_0,A_1,A_2>0$ are independent
of $\varepsilon$.

Next, we estimate the derivatives of
$\chi_\varepsilon$ on $S^{n-1}$.
Put
$$
h(\xi)=\left|
\frac{\xi}{|\xi|}-\widehat\xi_0
\right|^2.
$$
Using orthonormal coordinates for $|\cdot|$,
direct differentiation gives
$$
h_{\xi_i}
=2\left(
\langle\xi,\widehat\xi_0\rangle\xi_i
-(\widehat\xi_0)_i
\right)
$$
and
$$
\begin{aligned}
h_{\xi_i\xi_j}
=2\bigl(
(\widehat\xi_0)_i\xi_j
+(\widehat\xi_0)_j\xi_i+\langle\xi,\widehat\xi_0\rangle\delta_{ij}
-3\langle\xi,\widehat\xi_0\rangle
\xi_i\xi_j
\bigr).
\end{aligned}
$$
By the Cauchy--Schwarz inequality and the
support condition above, we have
$$
\begin{aligned}
\sum_{i=1}^n
\left(
\langle\xi,\widehat\xi_0\rangle\xi_i
-(\widehat\xi_0)_i
\right)^2
=1-\langle\xi,\widehat\xi_0\rangle^2
\leq|\xi-\widehat\xi_0|^2
\leq\varepsilon^2.
\end{aligned}
$$
Consequently,
$$
|h_{\xi_i}|\leq2\varepsilon,
\qquad
|h_{\xi_i\xi_j}|\leq12
$$
on $\operatorname{supp}
(\chi_\varepsilon|_{S^{n-1}})$.

Let $M_1,M_2>0$ be constants satisfying
$$
|\rho'(s)|\leq M_1,\qquad
|\rho''(s)|\leq M_2,\qquad s\geq0.
$$
Since
$\chi_\varepsilon=\rho(h/\varepsilon^2)$,
$$
(\chi_\varepsilon)_{\xi_i}
=\rho'\left(\frac{h}{\varepsilon^2}\right)
\frac{h_{\xi_i}}{\varepsilon^2},
$$
and
$$
\begin{aligned}
(\chi_\varepsilon)_{\xi_i\xi_j}
={}&
\rho''\left(\frac{h}{\varepsilon^2}\right)
\frac{h_{\xi_i}h_{\xi_j}}{\varepsilon^4}+
\rho'\left(\frac{h}{\varepsilon^2}\right)
\frac{h_{\xi_i\xi_j}}{\varepsilon^2}.
\end{aligned}
$$
Therefore,
$$
\begin{aligned}
\sup_{|\xi|=1}
|(\chi_\varepsilon)_{\xi_i}|
&\leq\frac{2M_1}{\varepsilon},\\
\sup_{|\xi|=1}
|(\chi_\varepsilon)_{\xi_i\xi_j}|
&\leq\frac{4M_2+12M_1}{\varepsilon^2}.
\end{aligned}
$$

Combining the preceding estimates,
we obtain
$$
\begin{aligned}
|\widetilde R_{\xi_i\xi_j}|
\leq{}&
|\chi_\varepsilon R_{\xi_i\xi_j}|
+|(\chi_\varepsilon)_{\xi_i}R_{\xi_j}|+|(\chi_\varepsilon)_{\xi_j}R_{\xi_i}|
+|R(\chi_\varepsilon)_{\xi_i\xi_j}|\\
\leq{}&
A_2\varepsilon
+4M_1A_1\varepsilon+A_0(4M_2+12M_1)\varepsilon\\
={}&C\varepsilon,
\end{aligned}
$$
where
$$
C=A_2+4M_1A_1
+A_0(4M_2+12M_1)
$$
is independent of $\varepsilon$.
Outside the support of $\chi_\varepsilon$,
the function $\widetilde R$ vanishes
locally. Hence,
\begin{equation}\label{eq:perturbation-bound}
\sup_{|\xi|=1}
|\widetilde R_{\xi_i\xi_j}(\xi)|
\leq C\varepsilon.
\end{equation}

Let $\lambda_*>0$ be the smallest
eigenvalue of the positive definite
matrix $\bigl(2g^{*ij}(\xi_0)\bigr)$.
By \eqref{eq:perturbation-bound},
for any $z=(z_1,\ldots,z_n)\in\mathbb{R}^n$
and $|\xi|=1$, we have
$$
\begin{aligned}
2\sum_{i,j=1}^n \widetilde{g}^{*ij}(\xi)z_i z_j
&=2\sum_{i,j=1}^n g^{*ij}(\xi_0)z_i z_j
+\sum_{i,j=1}^n\widetilde R_{\xi_i\xi_j}z_i z_j\\
&\geq \lambda_*|z|^2
-C\varepsilon\sum_{i,j=1}^n|z_i||z_j|\\
&\geq(\lambda_*-nC\varepsilon)|z|^2.
\end{aligned}
$$
Choosing $\varepsilon>0$ sufficiently small
so that $nC\varepsilon<\lambda_*/2$,
we obtain
$$
\sum_{i,j=1}^n \widetilde{g}^{*ij}(\xi)z_i z_j
\geq\frac{\lambda_*}{4}|z|^2>0
$$
for every $z\neq0$ and $|\xi|=1$. Since $\widetilde{g}^{*ij}(\xi)$ is zero-homogeneous, the
matrix $\bigl(\widetilde{g}^{*ij}(\xi)\bigr)$ is positive definite on
$\mathbb{R}^n\setminus\{0\}$.
Moreover,
$
\widetilde F^{*2}(\xi)
=\widetilde{g}^{*ij}(\xi)\xi_i\xi_j>0$ for $\xi\neq0$. Therefore, $\widetilde F^*$ is a regular
Minkowski norm on $\mathbb{R}^n$.

Finally, define
$$
C_0^*=
\left\{
\xi\in\mathbb{R}^n\setminus\{0\}:
\left|
\frac{\xi}{|\xi|}
-\frac{\xi_0}{|\xi_0|}
\right|<\frac{\varepsilon}{2}
\right\}.
$$
Then
$\mathbb{R}_{>0}\xi_0\subset C_0^*
\subset C_n^*$.
Since $\chi_\varepsilon=1$ on $C_0^*$,
we have
$$
\widetilde F^{*2}
=F_0^{*2}+\widetilde R
=F_0^{*2}+R
=F^{*2}
\quad\text{on }C_0^*.
$$
Hence $\widetilde F^*=F^*$ on $C_0^*$.
This completes the proof.
\endproof

Fix $(u_0,v_0)$ with $0<u_0<1/2$, and let
$\xi_0$ be a conormal vector to the helicoid
at this point, normalized so that
$F_{\mathrm{conic}}^*(\xi_0)=1$.
Put $\nu_0=(\xi_0,0)\in C_n^*$.
Applying Lemma \ref{lemX} to the
$n$-dimensional conic dual Minkowski norm $F^*$
at $\nu_0$, we obtain a regular dual Minkowski
norm $\widetilde F^*$ on $\mathbb R^n$ such that
$$
\widetilde F^*=F^*
\quad\text{on }C_0^*,
$$
where $C_0^*\subset C_n^*$ is an open cone
containing $\mathbb R_{>0}\nu_0$.
Let $\widetilde F$ denote the Minkowski norm
dual to $\widetilde F^*$.

Consider the cylinder over the helicoid,
$$
\widetilde{\mathbf r}(u,v,x_2)
=(u\cos v,u\sin v,av,x_2),
\qquad x_2\in\mathbb R^{n-3}.
$$
Its positive conormal field is given by
$$
\nu(u,v,x_2)
=(a\sin v,-a\cos v,u,0).
$$
Since
$\nu(u_0,v_0,0)\in\mathbb R_{>0}\nu_0
\subset C_0^*$,
the openness of $C_0^*$ and the continuity
of $\nu$ imply that
$$
\nu(u,v,x_2)\in C_0^*
$$
for $(u,v)$ in a sufficiently small
neighborhood of $(u_0,v_0)$.
Notice that $\nu$ is independent of $x_2$.

On $C_0^*$, we have $\widetilde{\mathbf n}=\mathcal L_{\widetilde F}^{-1}
\left(\frac{\nu}{\widetilde F^*(\nu)}\right)=\mathcal L_F^{-1}
\left(\frac{\nu}{F^*(\nu)}\right)=\mathbf n.$
Consequently, for every tangent vector $X$,
$$
\widetilde A(X)
=-\partial_X\widetilde{\mathbf n}
=-\partial_X\mathbf n
=A_{F}(X)=(A_{\mathrm{conic}}\oplus 0_{n-3})(X).
$$
Hence, by Lemma~\ref{thm17}, the principal curvatures are
$$
1,\qquad -1,\qquad 0,
$$
with multiplicities $1$, $1$, and $n-3$, respectively. Since all principal curvatures are constant,
Definition \ref{defi2} implies that this hypersurface is locally isoparametric. Moreover, their sum is zero, so it is minimal. Thus we obtain
\begin{theo}\label{theo42}
For every $n\geq4$, there exists a regular Minkowski metric on $\mathbf R^{n}$ admitting a local minimal isoparametric hypersurface with three distinct principal curvatures $1,-1,0$, of multiplicities $1,1,n-3$, respectively.
\end{theo}
\begin{rema}\label{rema1}
There exist local minimal isoparametric hypersurfaces in Minkowski spaces with $g=3$, which are not Minkowski hyperplanes, hyperspheres, or cylinders. Hence, unlike the Riemannian space form case, local isoparametric hypersurfaces in Finsler space forms need not admit global isoparametric extensions.
\end{rema}
\subsection{Umbilic hypersurfaces}
~~~~In \cite{HYS}, we know that any totally umbilic hypersurface in a Minkowski space must be isoparametric. In this section, we generalize the previous result to the following,
\begin{theo} \label{thm01}
Let $(N,F)$ be an $n$-dimensional Finsler space with constant flag curvature and $n\geq3$. Then every connected totally umbilic hypersurface is locally isoparametric.
\end{theo}
\proof
Suppose that $M$ is a totally umbilic hypersurface of $N$ with principal curvature $\lambda$, from (\ref{2.X}) and (\ref{2.5}),
\begin{align}\label{5.1}
D^{\textbf{n}}_{\partial u_{a}}\textbf{n}=-\lambda d\phi(\partial u_{a}).
\end{align}
Set $\textbf{n}=n^{i}\frac{\partial}{\partial x^{i}}$, $n^{i}_{a}=\frac{\partial n^{i}}{\partial u^{a}}$ and $n^{i}_{ab}=\frac{\partial n^{i}}{\partial u^{a}\partial u^{b}}$. From (\ref{new}),
\begin{align}\label{5.3}
D^{\textbf{n}}_{\partial u_{b}}D^{\textbf{n}}_{\partial u_{a}}\textbf{n}&=D^{\textbf{n}}_{\partial u_{b}}(n^{j}_{a}
+n^{i}\phi^{l}_{a}\Gamma^{j}_{li}(\textbf{n}))\partial x_{j}=\Big(n^{j}_{ab}+n^{i}\phi^{l}_{ab}\Gamma^{j}_{li}(\textbf{n})+(n^{i}_{b}\phi^{l}_{a}+n^{l}_{a}\phi^{i}_{b})\Gamma^{j}_{li}(\textbf{n})\nonumber\\
&+n^{i}\phi^{l}_{a}\phi^{p}_{b}\big(\frac{\partial\Gamma^{j}_{li}}{\partial x^{p}}(\textbf{n})+\frac{\partial\Gamma^{j}_{il}}{\partial y^{q}}(\textbf{n})\frac{\partial n^{q}}{\partial x^{p}}+\Gamma^{k}_{li}(\textbf{n})\Gamma^{j}_{pk}(\textbf{n})\big)\Big)\partial x_{j}.
\end{align}
From (\ref{2.6}), (\ref{2.6.6}) and (\ref{5.1}), we have
\begin{align}\label{5.3.3}
n^{i}\phi^{l}_{a}\phi^{p}_{b}(\frac{\partial\Gamma^{j}_{li}}{\partial x^{p}}(\textbf{n})+\frac{\partial\Gamma^{j}_{il}}{\partial y^{q}}(\textbf{n})\frac{\partial n^{q}}{\partial x^{p}})&=n^{i}\phi^{l}_{a}\phi^{p}_{b}\frac{\delta\Gamma^{j}_{li}}{\delta x^{p}}(\textbf{n})+n^{i}\phi^{l}_{a}\phi^{p}_{b}(N^{q}_{p}+\frac{\partial n^{q}}{\partial x^{p}})\frac{\partial \Gamma^{j}_{li}}{\partial y^{q}}(\textbf{n})\nonumber\\
&=n^{i}\phi^{l}_{a}\phi^{p}_{b}\frac{\delta\Gamma^{j}_{li}}{\delta x^{p}}(\textbf{n})+n^{i}\phi^{l}_{a}\phi^{p}_{b}(\nabla^{\textbf{n}}_{\partial x^{p}}\textbf{n})^{q}\frac{\partial \Gamma^{j}_{li}}{\partial y^{q}}(\textbf{n})\nonumber\\
&=n^{i}\phi^{l}_{a}\phi^{p}_{b}\frac{\delta\Gamma^{j}_{li}}{\delta x^{p}}(\textbf{n})+\lambda n^{i}P^{j}_{i~lq}(\textbf{n})\phi^{l}_{a}\phi^{q}_{b}\nonumber\\
&=n^{i}\phi^{l}_{a}\phi^{p}_{b}\frac{\delta\Gamma^{j}_{li}}{\delta x^{p}}(\textbf{n})-\lambda L^{j}_{~lq}(\textbf{n})\phi^{l}_{a}\phi^{q}_{b}.
\end{align}
Let $c$ be the constant flag curvature of $N$.
Then
$$
\begin{aligned}
R^j{}_p(y)=c\bigl(F^2(y)\delta^j_p-y^jy_p\bigr),~~~~y^iR^j{}_{ipl}(y)=c\bigl(\delta^j_p y_l-\delta^j_l y_p\bigr),
\end{aligned}
$$
where $y_i=g_{ij}(y)y^j$. Hence, from (\ref{2.6.6.6}), (\ref{5.3}) and (\ref{5.3.3}),
\begin{align}\label{5.5}
D^{\textbf{n}}_{\partial u_{b}}D^{\textbf{n}}_{\partial u_{a}}\textbf{n}-D^{\textbf{n}}_{\partial u_{a}}D^{\textbf{n}}_{\partial u_{b}}\textbf{n}
&=n^{i}\phi^{l}_{a}\phi^{p}_{b}R^{j}_{ipl}(\textbf{n})\partial x_{j}\nonumber\\
&=c\phi^{j}_{b}g_{\textbf{n}}(\textbf{n},d\phi(\partial u_{a}))-c\phi^{j}_{a}g_{\textbf{n}}(\textbf{n},d\phi(\partial u_{b}))\partial x_{j}=0.
\end{align}
From (\ref{5.1}) and (\ref{5.5}), by a direct calculation,
$$0=D^{\textbf{n}}_{\partial u_{b}}(-\lambda d\phi(\partial u_{a}))-D^{\textbf{n}}_{\partial u_{a}}(-\lambda d\phi(\partial u_{b}))=(\frac{\partial\lambda}{\partial u^{a}}\phi^{j}_{b}-\frac{\partial\lambda}{\partial u^{b}}\phi^{j}_{a})\partial x_{j}.$$
Since $n\geq3$, the tangent space of $M$ has dimension at least two. For each index $a$, choose $b\neq a$. Because $d\phi(\partial u_{a})$ and $d\phi(\partial u_{b})$ are linearly independent, the last identity gives $\partial\lambda/\partial u^{a}=\partial\lambda/\partial u^{b}=0$. Hence $d\lambda=0$ on $M$, and connectedness implies that $\lambda$ is a single constant on $M$. From Definition \ref{defi2}, $M$ is a local isoparametric hypersurface. This completes the proof.
\endproof
\section{Proof of Theorem \ref{thm03}}
Let $M_{0}=f^{-1}(t_{0})$ be the connected totally umbilic regular level in Theorem \ref{thm03}. By Theorem \ref{thm01}, its principal curvature is a constant $\lambda_{0}$. Since $f$ has no critical values in
$\operatorname{int}f(N)$, every maximal
$f$-segment intersects each regular level
exactly once and orthogonally (see Lemma 3.5 in \cite{CH}). Hence, the normal exponential map along these
$f$-segments induces a diffeomorphism between any two regular levels. Since $M_0$ is connected,
so is every regular level.

For $t_1<t_2$ in $\operatorname{int}f(N)$,
let $M_i=f^{-1}(t_i)$ $(i=1,2)$ and set
$$
s=\int_{t_1}^{t_2}
\frac{d\tau}{\sqrt{a(\tau)}}.
$$
Fix $x\in M_1$ and put
$y=\Phi_s(x)\in M_2$.
Let $A_1$ and $A_2$ be the shape operators
of $M_1$ at $x$ and $M_2$ at $y$,
respectively, with respect to
$\nabla f/F(\nabla f)$.
Denote by $P_s:T_xN\to T_yN$
the parallel translation along the
$f$-segment $\gamma_x$ from $x$ to $y$.

For the following calculations, define
\begin{equation}\label{eq11}
\mathfrak{s}_c(s)=
\begin{cases}
\dfrac{\sin(\sqrt{c}\,s)}{\sqrt{c}},
& c>0,\\[6pt]
s, & c=0,\\[6pt]
\dfrac{\sinh(\sqrt{-c}\,s)}{\sqrt{-c}},
& c<0,
\end{cases}
\qquad
\mathfrak{c}_c(s)=\mathfrak{s}_c'(s).
\end{equation}

For $X\in T_xM_1$, let $J$ be the Jacobi field induced by the normal variation $\Phi_{s}(x(t))$ as in the proof of Theorem \ref{thm3}. Then
$$
J(0)=X,\qquad D_T^T J(0)=-A_1X.
$$
Since $J$ is orthogonal to $T$ and the
flag curvature is constant $c$, we have
$$
D_T^T D_T^T J+cJ=0.
$$
Setting $Y(s)=P_s^{-1}J(s)$, we obtain
$$
Y''(s)+cY(s)=0,\qquad
Y(0)=X,\quad Y'(0)=-A_1X.
$$
Hence
\begin{equation}\label{eqq}
J(s)=P_s\left(
\mathfrak{c}_c(s)I-\mathfrak{s}_c(s)A_1
\right)X.
\end{equation}
Consequently,
$$
(d\Phi_s)_x
=P_s\left(
\mathfrak{c}_c(s)I-\mathfrak{s}_c(s)A_1
\right).
$$

Since $\Phi_s:M_1\to M_2$ is a diffeomorphism, $(d\Phi_s)_x$ is a
linear isomorphism. As $P_s$ is also
a linear isomorphism,
$$
B_s:=\mathfrak{c}_c(s)I
-\mathfrak{s}_c(s)A_1
$$
is invertible. Thus
$$
T_yM_2=(d\Phi_s)_x(T_xM_1)=P_sB_s(T_xM_1)=P_s(T_xM_1).
$$
So $P_s$ restricts to a linear
isomorphism from $T_xM_1$ onto $T_yM_2$.

Differentiating (5.2), we obtain
$$
-D_T^T J(s)
=P_s\left(
c\mathfrak{s}_c(s)I+
\mathfrak{c}_c(s)A_1
\right)X.
$$
Since $A_2J(s)=-D_T^TJ(s)$ and
$B_s$ is invertible, we have
\begin{equation}\label{eq1}
P_s^{-1}A_2P_s=
\left(
c\mathfrak{s}_c(s)I+
\mathfrak{c}_c(s)A_1
\right)
\left(
\mathfrak{c}_c(s)I-
\mathfrak{s}_c(s)A_1
\right)^{-1}.
\end{equation}
Consequently, if $\kappa_1$ is a principal
curvature of $M_1$, the corresponding
principal curvature of $M_2$ is
\begin{equation}\label{eq2}
\kappa_2(s)=
\frac{
c\mathfrak{s}_c(s)+
\mathfrak{c}_c(s)\kappa_1
}{
\mathfrak{c}_c(s)-
\mathfrak{s}_c(s)\kappa_1
}.
\end{equation}

By (\ref{eq2}), we obtain
$$
\kappa(t)=
\begin{cases}
\displaystyle
\frac{\lambda_0\mathfrak{c}_c(s)
+c\mathfrak{s}_c(s)}
{\mathfrak{c}_c(s)
-\lambda_0\mathfrak{s}_c(s)},
& t>t_0,\\[12pt]
\displaystyle
\frac{\lambda_0\mathfrak{c}_c(s)
-c\mathfrak{s}_c(s)}
{\mathfrak{c}_c(s)
+\lambda_0\mathfrak{s}_c(s)},
& t<t_0,
\end{cases}
$$
where
$$
s=\int_{\min\{t,t_0\}}^{\max\{t,t_0\}}
\frac{d\tau}{\sqrt{a(\tau)}}.
$$
Thus every regular level is connected and
totally umbilic with constant principal
curvature. By Definition~2.2, $f$ is
isoparametric. This completes the proof.
\begin{rema}\label{rema2}
In contrast to Remark \ref{rema1}, Theorems \ref{thm03} and \ref{thm01} show that, if $f$ is a transnormal function on $(\overrightarrow{N}(c),F)$, then the local isoparametric hypersurfaces can be extended to the global isoparametric hypersurfaces under additional assumptions.
\end{rema}
\section{Proof of Theorem \ref{thm02}}
~~~~From \cite{HYS}, it is easy to get that
\begin{lemm}\label{lem4-8} Let $f$ be an isoparametric function on a Finsler manifold $(N,F)$ with constant flag curvature $c$. Then the principal curvature $\kappa_{i}$ of an isoparametric hypersurface $f^{-1}(t)$ satisfying
 \begin{align}
\sum_{i=1}^{n-1}\kappa_i&=\frac{a'(t)}{2\sqrt{a(t)}}-\frac{b(t)}{\sqrt{a(t)}},\label{3.41}\\
\sqrt{a(t)}\frac{\partial \kappa_i}{\partial t}&=c+\kappa_i^2,\quad
i=1,\cdots,n-1.\label{3.81}\end{align}
\end{lemm}

Let~$r(x)= d_F(p, x)$ be the distance function from $p$ to $x$ on $(\overrightarrow{N}(c),F)$. Set $f=\frac{1}{2}r^{2}$.
Obviously, $f$ is a non-constant $C^1$ function on $\mathcal{D}_{p}$, which is smooth on $\mathcal{D}_{p}\setminus\{p\}$. In particular, according to Cartan-Hadamard Theorem, $f$ is a globally defined $C^1$ function on $(\overrightarrow{N}(c),F)$ for $c \leq 0$, retaining its smoothness away from the point $p$.

By a direct calculation, $F(\nabla f)=\sqrt{2f}$. From the proof of Theorem 4.1 in \cite{WX}, we have $\hat{\Delta} r=(n-1)\frac{\mathfrak{s}'_{c}(r)}{\mathfrak{s}_{c}(r)}$, where $\mathfrak{s}_{c}$ is given by (\ref{eq11}). Hence,
$$\hat{\Delta} f=r\hat{\Delta}r+1=r(n-1)\frac{\mathfrak{s}'_{c}(r)}{\mathfrak{s}_{c}(r)}+1=\sqrt{2f}(n-1)\frac{\mathfrak{s}'_{c}(\sqrt{2f})}{\mathfrak{s}_{c}(\sqrt{2f})}+1.$$
Hence, $f$ satisfies (\ref{1.1}) with $a(t)=2t, b(t)=\sqrt{2t}(n-1)\frac{\mathfrak{s}'_{c}(\sqrt{2t})}{\mathfrak{s}_{c}(\sqrt{2t})}+1$. Therefore, $f$ is a local isoparametric function on $\mathcal{D}_{p}$, and when $c\leq0$, it is a global isoparametric function on $(\overrightarrow{N}(c),F)$.

By (\ref{3.41}) and (\ref{3.81}), a
straightforward calculation yields that
\begin{align}\label{4.6}
\sum\limits_{a}\kappa_{a}=\left\{
\begin{array}{lcl}
\frac{1}{\sqrt{2t}}(1-n),       &      & {c=0,}\\
-(n-1)\sqrt{c}\cot(\sqrt{2ct}),    &      & {c>0,}\\
-(n-1)\sqrt{-c}\coth(\sqrt{-2ct}),       &      & {c<0,}
\end{array} \right.
\end{align}
and
\begin{align}\label{4.7}
\sum\limits_{a}\kappa_{a}^2=\left\{
\begin{array}{lcl}
\frac{1}{2t}(n-1),       &      & {c=0,}\\
(n-1)c(\csc^2(\sqrt{2ct})-1),    &      & {c>0,}\\
-c(n-1)(\textmd{csch}^2(\sqrt{-2ct})+1),       &      & {c<0.}
\end{array} \right.
\end{align}
(\ref{4.6}) and (\ref{4.7}) imply that
  \begin{align}\label{5.32}
(\sum_a \kappa_a)^2&=(n-1)\sum_a \kappa^2_a.
\end{align}
From (\ref{4.6}) and (\ref{5.32}), we have
\begin{align*}
\kappa_1=\kappa_2=\ldots=\kappa_{n-1}=\left\{
\begin{array}{lcl}
-\frac{1}{\sqrt{2t}},       &      & {c=0,}\\
-\sqrt{c}\cot(\sqrt{2ct}),    &      & {c>0,}\\
-\sqrt{-c}\coth(\sqrt{-2ct}),       &      & {c<0.}
\end{array} \right.
\end{align*}
That is, the principal curvature $\kappa_a(t)$ is a constant for $M_t \cap \mathcal{D}_p$, and for $M_t$ if $c \leq 0$.

Let $\mathbb{S}^{n}$ be an $n$-sphere equipped with a Finsler metric $F$ and the Busemann-Hausdorff volume form $d\mu$. If it has constant flag curvature $c$, vanishing $\textbf{S}$-curvature and $Diam=\frac{\pi}{\sqrt{c}}$, we call it a \textit{standard Finsler sphere} and denote it by $(\mathbb{S}^{n}(c),F,d\mu)$\cite{YH}. We can obtain a global isoparametric function on $(\mathbb{S}^{n}(c),F,d\mu)$ with $g=1$. Set $f(x)=-\cos(\sqrt{c}r(x))$. From \cite{YH}, we know that $f(x)$ is the first eigenfunction of $(\mathbb{S}^{n}(c),F,d\mu)$. Hence, from \cite{GS}, $f\in C^{\infty}(\mathbb{S}^{n}_{f})$ and $f\in C^{1}(\mathbb{S}^{n})$. In addition, since $F(\nabla f)=\sqrt{c(1-f^2)}$ and $\Delta f=-ncf$, then $f$ is isoparametric on $(\mathbb{S}^{n}(c),F,d\mu)$ with $a(t)=c(1-t^2)$ and $b(t)=-nct$. From Lemma \ref{lem4-8}, we also have $\kappa_{1}=\kappa_{2}=\cdots=\kappa_{n-1}=\frac{\sqrt{c}t}{\sqrt{1-t^2}}$. That completes the proof of Theorem \ref{thm02}.

\section*{Declarations}

This work was supported by NNSFC(No. 12401056) and NNSFC(No. 12261034). The authors declare no conflict of interest. All authors contributed to the proof and calculation. The distinction between local isoparametric hypersurfaces and global isoparametric hypersurfaces in Finsler manifolds is proposed by Qun He. The first draft of the manuscript was written by Yali Chen and all authors commented on previous versions of the manuscript. All authors read and approved the final manuscript.

Yali Chen\\
School of Mathematics and Statistics, Anhui Normal University, Wuhu, Anhui, 241000, China\\
E-mail: chenyl@ahnu.edu.cn\\

Qun He \\
School of Mathematical Sciences, Tongji University, Shanghai, 200092, China\\
E-mail: hequn@tongji.edu.cn


\begin{thebibliography}{99}
\bibitem{GTM}
D. Bao, S. S. Chern, Z. Shen. An Introduction to Riemann-Finsler Geometry. Graduate Texts in Mathematics 200. New York: Springer-Verlag, 2000.
\bibitem{C}
E. Cartan. Familles de surfaces isoparam$\acute{e}$triques dans les espaces $\grave{a}$ courbure constante.  Ann. Mat. Pura Appl., 1938, 17: 177-191.
\bibitem{C1}
Q. S. Chi. Isoparametric hypersurfaces with four principal curvatures, IV. J. Differential Geom., 2020, 115(2): 225-301.
\bibitem {TP}
T. E. Cecil and P. J. Ryan. Geometry of Hypersurfaces. Springer Monographs in Math., 2015.
\bibitem{CS1}
T. E. Cecil, Q. S. Chi and G. R. Jensen. Isoparametric hypersurfaces with four principal curvatures. Ann. Math, 2007, 166: 1-76.
\bibitem{CH}
Y. L. Chen and Q. He. Transnormal Functions and Focal Varieties on Finsler Manifolds. J. Geom. Anal., 2023, 33(4): 128.
\bibitem{HD}
P. L. Dong and Q. He. Isoparametric hypersurfaces of a class of Finsler manifolds induced by navigation problem in Minkowski spaces. Diff. Geom. Appl., 2020, 68: 101581.
\bibitem{DC}
P. L. Dong and Y. L. Chen. Isoparametric hypersurfaces and hypersurfaces with constant principal curvatures in Finsler spaces. Ann. Pol. Math., 2023, 131(2): 127-140.
\bibitem {GS} Y. Ge and Z. Shen. Eigenvalues and eigenfunctions of metric measure manifolds. Proc. London Math., 2001, 82(3): 725-746.
\bibitem{HYS}
Q. He, S. T. Yin and Y. B. Shen. Isoparametric hypersurfaces in Minkowski spaces. Diff. Geom. Appl., 2016, 47: 133-158.
\bibitem{HYS1}
Q. He, S. T. Yin and Y. B. Shen. Isoparametric hypersurfaces in Funk manifolds. Science China Math., 2017, 60(12): 2447-2464.
\bibitem{HDY}
Q. He, P. L. Dong and S. T. Yin. Isoparametric hypersurfaces in Randers space forms. Acta Math. Sin. Engl. Ser., 2020, 36(9): 1049-1060.
\bibitem{HCY}
Q. He, Y. L. Chen, S. T. Yin and T. T. Ren. Isoparametric hypersurfaces in Finsler space forms. Science China Math., 2021, 64(7): 1463-1478.
\bibitem{GM}
J. Q. Ge and H. Ma. Anisotropic isoparametric hypersurfaces in Euclidean spaces. Ann. Global Anal. Geom., 2012, 41: 347--355.
\bibitem{HH}
Q. He, X. Huang and P. L. Dong. Isoparametric hypersurfaces in conic Finsler manifolds. Differ. Geom. Appl., 2022, 84: 101937.
\bibitem{M}
R. Miyaoka. Transnormal functions on a Riemannian manifold. Diff. Geom. Appl., 2013, 31(1): 130-139.
\bibitem{QM}
Q. M. Wang. Isoparametric functions on Riemannian manifolds I. Math. Ann., 1987, 277: 639-646.
\bibitem {WX} B. Y. Wu and Y. L. Xin. Comparison theorems in Finsler geometry
and their applications. Math. Ann., 2007, 337: 177-196.
\bibitem{SZ}
Z. M. Shen. Lectures on Finsler geometry. Singapore: World Scientific Publishing Co 2001.
\bibitem{SZ1}
Z. M. Shen. On Finsler geometry of submanifolds. Math. Ann., 1998, 311(3): 549-576.
\bibitem{YH}
S. T. Yin and Q. He. The maximum diam theorem on Finsler manifolds. J. Geom. Anal., 2021, 31: 12231-12249.
\end{thebibliography}
\end{document}